\documentclass[11pt]{article}
\usepackage{amscd, amsthm, color}
\usepackage{mathrsfs}
\input epsf
\usepackage{amssymb}
\usepackage{amsmath}
\usepackage{amsfonts}
\input amssym.def
\newtheorem{theorem}{Theorem}[section]

\newtheorem{lemma}[theorem]{Lemma}
\newtheorem{prop}[theorem]{Proposition}

\newtheorem{re}[theorem]{Remark}
\newtheorem{no}[theorem]{Notation}
\newtheorem{definition}[theorem]{Definition}
\theoremstyle{definition}

\definecolor{wco}{rgb}{0.5,0.2,0.3}

\numberwithin{equation}{section}

\begin{document}
\title{Gelfand and Kolmogorov numbers of Sobolev embeddings\\
of weighted function spaces
}
\author{Shun Zhang $^{a,\, b}$,\quad \ Gensun Fang
$^{a,\,}$\footnote{Corresponding author.
\newline\indent\ \, E-mail addresses: fanggs@bnu.edu.cn (G. Fang), shzhang27@163.com (S. Zhang).}
\\ {\small $^{a}$ School of Mathematical Sciences, Beijing Normal University,
Beijing 100875, China}
\\ {\small $^{b}$ School of Computer Science and Technology, Anhui
University,
 Hefei 230039, Anhui, China}}
\maketitle {\bf Abstract.}\, In this paper we study the Gelfand and
Kolmogorov numbers of Sobolev embeddings between weighted function
spaces of Besov and Triebel-Lizorkin type with polynomial weights.
The sharp asymptotic estimates are determined in
 the so-called non-limiting case. 

\section{Introduction}
~~~~ In recent years a great deal of effort has gone into studying
compactness of Sobolev embeddings between function spaces of Besov
and Triebel-Lizorkin type from the standpoint of n-widths,
especially approximation, Gelfand and Kolmogorov numbers. The case
of function spaces defined on bounded domains has attracted a lot of
attention, see \cite{ET96,Pin85,Sk06,Tr06,Vy08}. For weighted
function spaces of this type, such embeddings have also been studied
by many authors, with polynomial (and more general) weights
considered. Some breakthroughs on approximation numbers may be found
in the works of Caetano \cite{Ca}, Haroske \cite{Ha95,Ha97}, Haroske
and Skrzypczak \cite{HS08,HS11}, Mynbaev and Otel'baev \cite{MO88},
Piotrowska \cite{Pio08}, Skrzypczak \cite{Sk05,Sk06} and Vasil'eva
\cite{Va10}. In particular, Skrzypczak \cite{Sk05} investigated the
approximation numbers of the embeddings in the case of polynomial
weights, by using operator ideals. In the context of Gelfand and
Kolmogorov numbers, Vasil'eva \cite{Va08,Va10} established  the
asymptotics of the Kolmogorov numbers of weighted Sobolev classes on
a finite interval or half-axis in the space $L_q$ with weight.
However, the estimates in many other cases are still left open.

$n$-Widths are a well-explored subject in approximation theory, see
\cite{Pie78,Pin85,Vy08}, and recently they have been applied in many
areas, including compressed sensing \cite{Do}, computational
mechanics \cite{BBO,EBBH} and spectral theory \cite{ET96}. In
particular, in the remarkable paper introducing compressed sensing
\cite{Do}, general performance bounds for sparse recovery methods
are obtained by means of the theory of $n$-widths. In \cite{BBO,EBBH}
the Kolmogorov number is utilized to assess approximation properties
of functions employed in finite element techniques. The performance
of approximation numbers for describing the spectral properties of
(pseudo-)differential operators is discussed in \cite{ET96}.

In this paper we present the sharp asymptotic estimates of the
Gelfand and Kolmogorov numbers in the so-called non-limiting case.
 Although there are parallel considerations which cover some cases in
Vasil'eva \cite{Va08,Va10}, we proceed in a completely different way
than Vasil'eva. We shall follow the method utilized for
approximation numbers in Skrzypczak \cite{Sk05} with its corrigendum
\cite{SV09}.

Motivated by \cite{KLSS06b,Sk05,Sk06}, using the discretization
method due to Maiorov \cite{Mai75}, we reduce the function space of
the problem to a weighted sequence space, then we determine the
asymptotic behavior of the Gelfand and Kolmogorov numbers of
Sobolev embeddings between weighted function spaces. The
discretization technique is very important in the process of
determining the exact order of n-widths of such classes, and in many
cases it plays a major role. Moreover, our main tools are the use of
operator ideals, see \cite{Car81,Pie78,Pie87}, and the basic
estimates of related widths of the Euclidean ball due to Gluskin
\cite{Gl83}. Historically, the technique of estimating single
n-widths via estimates of ideal quasi-norms derives from ideas of
Carl \cite{Car81}.

Following Skrzypczak \cite{Sk05}, we concentrate on the spaces with
polynomial weights
\begin{equation}\label{w_a}
w_\alpha(x):=(1+|x|^2)^{\alpha/2}
\end{equation}
for some exponent $\alpha>0.$ Let
\begin{equation}\label{emB_con} -\infty<s_2<s_1<\infty,\quad\quad
1\leq p_1\leq p_2\leq\infty, \quad\quad 1\leq q_1, q_2\leq\infty.
\end{equation}
It is well known that if $\delta=s_1-s_2-d(\frac 1{p_1}-\frac
1{p_2})>0$ then
\begin{equation}\label{emB_wei}
B_{p_1,q_1}^{s_1}(\mathbb{R}^d, w_\alpha)\hookrightarrow
B_{p_2,q_2}^{s_2}(\mathbb{R}^d).
\end{equation}
Moreover, the non-limiting case means as usual that
$\delta\neq\alpha$. The case $p_2< p_1$ is also considered in this
article.
\begin{re}\label{w12}
In \cite{HT94}, it is proved that the embedding
\begin{equation}\label{emB_wei12}
B_{p_1,q_1}^{s_1}(\mathbb{R}^d, v_1)\hookrightarrow
B_{p_2,q_2}^{s_2}(\mathbb{R}^d, v_2)
\end{equation}
(and its F-counterparts with $p_2 < \infty$) is compact if, and only
if,
\begin{equation}
s_1-\frac d{p_1}>s_2-\frac d{p_2}\quad{\rm and}\quad
\frac{v_2(x)}{v_1(x)}\rightarrow 0 \quad{\rm for}\quad
|x|\rightarrow \infty,
\end{equation}
where $-\infty<s_2<s_1<\infty, 0< p_1\le p_2\le \infty, 0< q_1,
q_2\leq\infty$ and $v_1, v_2$ are admissible weight functions, see
also \cite{Ha95}.

Based on these considerations (where comparatively general weight
functions are involved), we can assume that the target space is an
unweighted space. Therefore, we restrict ourselves to the so-called
standard situation:
$v_1(x)=w_\alpha(x)=(1+|x|^2)^{\alpha/2},\,\alpha>0,\,v_2(x)\equiv
1.$
\end{re}
\begin{no}
By the symbol ` $\hookrightarrow$'  we denote continuous embeddings.

Identity operators will always be denoted by {\rm id}. Sometimes we
do not indicate the spaces where {\rm id} is considered, and
likewise for other operators.

Let $X$ and $Y$ be complex Banach spaces and denote by $\mathcal
{L}(X, Y)$ the class of all linear continuous operators $T:\,X
\rightarrow\, Y.$ If no ambiguity arises, we write $\|T\|$ instead
of the more exact versions $\|T ~|~ \mathcal {L}(X, Y)\|$ or
$\|T:X\rightarrow Y\|$.

The symbol $a_n\sim b_n$ means that there exists a constant $c
> 0$ independent of $n$ such
that $$c^{-1}a_n\leq b_n\leq c a_n,\quad\quad n=1, 2, 3\ldots.$$

All unimportant constants will be denoted by $c$ or $C$, sometimes
with additional indices.
\end{no}

We start with recalling the definitions of Kolmogorov and Gelfand
numbers, cf. \cite{Pin85}. We use the symbol $A\subset\subset B$ if
$A$ is a closed subspace of a topological vector space $B$.
\begin{definition}
 Let $T\in\mathcal {L}(X,Y)$.
\begin{enumerate}
\item[$(i)$]\ The {\it $n$th Kolmogorov number} of the operator
$T$ is  defined by
\begin{equation*}
d_n(T, X, Y)=\inf\{\|Q_N^YT\|:\,N\subset\subset Y,\,{\rm dim}
(N)<n\},
\end{equation*}
also written by $d_n(T)$ if no confusion is possible. Here, $Q_N^Y$
stands for the natural surjection of\,\,\,$Y$ onto the quotient
space $Y/N$.

\item[$(ii)$]\ The {\it $n$th Gelfand number} of the operator
$T$ is  defined by
\begin{equation*}
c_n(T, X, Y)=\inf\{\|TJ_M^X\|:\,M\subset\subset X,\,{\rm codim}
(M)<n\},
\end{equation*}
also written by $c_n(T)$ if no confusion is possible. Here, $J_M^X$
stands for the natural injection of\,\,\,$M$ into $X$.
\end{enumerate}
\end{definition}

It is well-known that the operator $T$ is compact if and only if
$\lim_n d_n(T)=0$ or equivalently $\lim_n c_n(T)=0$, see
\cite{Pin85}.

The Kolmogorov and Gelfand numbers are both examples of so-called
$s$-numbers, cf. \cite{Pie78, Pie87, Pin85}. Let $s_n$ denote either
of these two quantities, $c_n$ or $d_n$, and let $Y$ be a Banach
space. We collect several common properties of Kolmogorov and
Gelfand numbers below: \vspace{-0.2cm}
\begin{enumerate}
\item[]{\rm\bf(PS1)}\ (nonincreasing property)\,
$\|T\|=s_1(T)\ge s_2(T)\ge\cdots\ge 0$ for all $T\in\mathcal{L}(X,
Y)$,\vspace{-0.2cm}

\item[]{\rm\bf(PS2)}\ (subadditivity)\, $s_{m+n-1}(S+T)\leq
s_m(S)+s_n(T)$\, for all $m, n\in\mathbb{N},\,\,S,
T\in\mathcal{L}(X, Y)$,\vspace{-0.2cm}

\item[]{\rm\bf(PS3)}\ (multiplicativity)\, $s_{m+n-1}(ST)\leq
s_m(S)s_n(T)$\, for all $T\in\mathcal{L}(X, Y)$, $S\in\mathcal{L}(Y,
Z)$

\quad\quad and $m, n\in\mathbb{N},$\, cf. \cite[p. 155]{Pie78},
where $Z$ denotes a Banach space,\vspace{-0.2cm}

\item[]{\rm\bf(PS4)}\ (rank property)\, ${\rm rank}(T)<n$ if and only if
$s_n(T)=0$, where $T\in\mathcal{L}(X, Y)$.\vspace{-0.2cm}
\end{enumerate}
Both concepts, Kolmogorov and Gelfand numbers, are related to each
other. Namely they are dual to each other in the following sense,
cf. \cite{Pie78, Pin85}: If $X$ and $Y$ are Banach spaces, then
\begin{equation}\label{dualc*d}
c_n(T^\ast)=d_n(T)
\end{equation}
for all compact operators $T\in\mathcal{L}(X, Y)$ and
\begin{equation}\label{duald*c}
d_n(T^\ast)=c_n(T)
\end{equation}
for all $T\in\mathcal{L}(X, Y).$

Following Pietsch \cite{Pie87}, we associate to the sequence of the
 Kolmogorov (or Gelfand) numbers the following operator ideals, and
for $0<r<\infty$, we put
\begin{equation}\mathscr{L}_{r,\infty}^{(s)}:=\left\{T\in\mathcal{L}(X,
Y):\quad \sup\limits_{n\in\mathbb{N}}n^{1/r}s_n(T)<\infty\right\}.
\end{equation}
 Equipped with the quasi-norm
\begin{equation}\label{idealddef}
L_{r,\infty}^{(s)}(T):=\sup\limits_{n\in\mathbb{N}}n^{1/r}s_n(T),
\end{equation}
the set $\mathscr{L}_{r,\infty}^{(s)}$ becomes a quasi-Banach space.
For such quasi-Banach spaces there always exists a real number
$0<\rho\leq 1$ such that
\begin{equation}\label{idealdinq}
L_{r,\infty}^{(s)}\left(\sum\limits_jT_j\right)^\rho\leq
C \sum\limits_jL_{r,\infty}^{(s)}(T_j)^\rho
\end{equation}
holds for any sequence of operators
$T_j\in\mathscr{L}_{r,\infty}^{(s)}.$ Then we shall use the
quasi-norms $L_{r,\infty}^{(c)}$ and $L_{r,\infty}^{(d)}$ for the
Gelfand and Kolmogorov numbers, respectively.

The paper is structured as follows. In Sect.~2, we introduce
weighted function spaces of $B$-type and $F$-type, and provide our main
results. In Sect.~3, the crucial part of the work will be done, we
investigate the Kolmogorov numbers of embeddings of related sequence
spaces. Finally, in Sect.~4, these results will be used to derive
the desired Kolmogorov number estimates for the function space
embeddings under consideration, and similar results on the Gelfand
numbers of such embeddings are established. Our main assertions are
Theorem \ref{kn} and Theorem \ref{gn}.
\section{Main results}
~~~~ We suppose that the reader is familiar with (unweighted)
function spaces of $B$-type and $F$-type on $\mathbb{R}^d$ . One can
consult \cite{ET96,Tr83} and many other literatures for the
definitions and basic properties.

Throughout this paper we are interested in the function spaces with
polynomial weights given by (\ref{w_a}). As usual,
$\mathcal{S}^\prime(\mathbb{R}^d)$ denotes the set of all tempered
distributions on the Euclidean $d$-space $\mathbb{R}^d$. For us it
will be convenient to introduce weighted function spaces to be
studied here.
\begin{definition}\label{BF}
Let $1\leq p\leq \infty,\,1\leq q\leq \infty,\,$ and
$s\in\mathbb{R}$. Then we put
\begin{equation*}
B_{p,q}^{s}(\mathbb{R}^d, w_\alpha)=\left\{f\in \mathcal{S}^\prime
(\mathbb{R}^d)\,:\, \|f ~|~ B_{p,q}^{s}(\mathbb{R}^d,
w_\alpha)\|=\|fw_\alpha ~|~
B_{p,q}^{s}(\mathbb{R}^d)\|<\infty\right\},
\end{equation*}
\begin{equation*}
F_{p,q}^{s}(\mathbb{R}^d, w_\alpha)=\left\{f\in \mathcal{S}^\prime
(\mathbb{R}^d)\,:\, \|f ~|~ F_{p,q}^{s}(\mathbb{R}^d,
w_\alpha)\|=\|fw_\alpha ~|~
F_{p,q}^{s}(\mathbb{R}^d)\|<\infty\right\},
\end{equation*}
with $p<\infty$ for the $F$-spaces.
\end{definition}
\begin{re}
If no ambiguity arises, then we can write $B_{p,q}^{s}(w_\alpha)$
and $F_{p,q}^{s}(w_\alpha)$ for brevity.
\end{re}
\begin{re}
There are different ways to introduce weighted function spaces, see,
e.g., Edmunds and Triebel \cite{ET96}, or Schmeisser and Triebel
\cite{ST87}. One can also consult \cite{ KLSS06a,Tr06} for related
remarks.
\end{re}

Let
$A_{p_,q}^{s}(\mathbb{R}^d,w_\alpha)\,(A_{p_,q}^{s}(\mathbb{R}^d))$
stand for either
$B_{p_,q}^{s}(\mathbb{R}^d,w_\alpha)\,(B_{p_,q}^{s}(\mathbb{R}^d))$
or $F_{p_,q}^{s}(\mathbb{R}^d,w_\alpha)$ \,
$(F_{p_,q}^{s}(\mathbb{R}^d)),$ with the constraint that for the
F-spaces $p<\infty$ holds.

Now we give a necessary and sufficient condition for compactness of
the embeddings under consideration, which was proved in \cite{HT94},
cf. also \cite{ET96,KLSS06b}.
\begin{prop}\label{compact}
Suppose $1\le p_1, p_2, q_1, q_2\le\infty\,$ and
$-\infty<s_2<s_1<\infty.$\, Let $\alpha>0,\,\delta=s_1-s_2-d(\frac
1{p_1}-\frac 1{p_2}).$ The embedding
$B_{p_1,q_1}^{s_1}(\mathbb{R}^d,w_\alpha)\hookrightarrow
B_{p_2,q_2}^{s_2}(\mathbb{R}^d)$ is compact if and only if
$\min(\alpha, \delta)>d \max(\frac 1{p_2}-\frac 1{p_1}, 0)$.
\end{prop}

A similar theorem also holds for $F_{p,q}^s$-spaces. We are now
ready to formulate our main results.
\begin{theorem}\label{kn}
Suppose $1\le p_1, p_2, q_1, q_2\le\infty\,$ and
$-\infty<s_2<s_1<\infty.$\, Let $\alpha>0,\,\delta=s_1-s_2-d(\frac
1{p_1}-\frac 1{p_2})>0,\,\theta =
\frac{1/{p_1}-1/{p_2}}{1/2-1/{p_2}},$\, and $\frac
1{\tilde{p}}=\frac\mu d+\frac 1{p_1},$ where
$\mu=\min(\alpha,\delta).\,$ Besides, we assume that\\
{\rm(a)}~\,$1\le p_1\le p_2\le\infty\,~or~\,\tilde{p}<p_2<p_1\leq\infty,$\\
{\rm(b)}~\,$\delta\neq\alpha,$\\
{\rm(c)}~\,$p_2<\infty$ ~when~ $p_1< p_2.$

Denote by $d_n$ the $n$th Kolmogorov number of the Sobolev embedding
\begin{equation}\label{aa}
A_{p_1,q_1}^{s_1}(\mathbb{R}^d,w_\alpha)\hookrightarrow
A_{p_2,q_2}^{s_2}(\mathbb{R}^d).
\end{equation}
Then

$\indent d_{n}\sim n^{-\varkappa},$\\ where\vspace{-0.2cm}
\begin{enumerate}
\item[{\rm (i)}]\ $\varkappa=\frac\mu d$\, if\, $1\le p_1\le
p_2\le 2$\,\,or\,\,$2<p_1 = p_2\le \infty$, \vspace{-0.2cm}

\item[{\rm (ii)}]\ $\varkappa =\frac \mu d+\frac 1{p_1}-\frac
1{p_2}$\, if\, $\tilde{p}<p_2<p_1\leq\infty$, \vspace{-0.2cm}

\item[{\rm (iii)}]\ $\varkappa=\frac \mu d +\frac 12-\frac 1{p_2}$\,
 if\, $1\le p_1 < 2 < p_2< \infty$\ and \,$\mu>\frac d{p_2}$,\vspace{-0.2cm}

\item[{\rm (iv)}]\ $\varkappa=\frac \mu d\cdot\frac{p_2}2$\,
 if\, $1\le p_1 < 2 < p_2< \infty$\ and \,$\mu<\frac d{p_2}$,\vspace{-0.2cm}

\item[{\rm (v)}]\ $\varkappa=\frac \mu d +\frac 1{p_1}-\frac 1{p_2}$\,
 if\, $2 \le p_1 < p_2 < \infty$\ and \,$\mu>\frac d{p_2}\theta$,\vspace{-0.2cm}

\item[{\rm (vi)}]\ $\varkappa=\frac \mu d\cdot\frac{p_2}2$\,
 if\, $2 \le p_1 < p_2 <\infty$\ and \,$\mu<\frac d{p_2}\theta$.
\end{enumerate}
\end{theorem}
\begin{re}
Similar conclusions on the $n$th Kolmogorov number could be made for
Corollary 19 in Skrzypczak \cite{Sk05} with its corrigendum to part
(iv) given in \cite{SV09}. Of course, the counterexample to our new
part (iv) could be also made for the limiting case $\delta=\frac
d{p_2}$ by virtue of the special example appeared at the end of
\cite{SV09}.
\end{re}

For $1\le p\le\infty,$ we set

$p^\prime=
\begin{cases}
\frac p{p-1}\quad &{\rm if}\ 1<p<\infty,\\
1 &{\rm if}\ p=\infty,\\
\infty &{\rm if}\ p=1.
\end{cases}
$

\begin{theorem}\label{gn}
Suppose $1\le p_1, p_2, q_1, q_2\le\infty\,$ and
$-\infty<s_2<s_1<\infty.$\, Let $\alpha>0,\,\delta=s_1-s_2-d(\frac
1{p_1}-\frac 1{p_2})>0,\,\theta_1 =
\frac{1/{p_1}-1/{p_2}}{1/{p_1}-1/2}$,\, and $\frac
1{\tilde{p}}=\frac\mu d+\frac 1{p_1},$ where
$\mu=\min(\alpha,\delta).\,$ Besides, we assume that\\
{\rm(a)}~\,$1\le p_1\le p_2\le\infty\,~or~\,\tilde{p}<p_2<p_1\leq\infty,$\\
{\rm(b)}~\,$\delta\neq\alpha,$\\
{\rm(c)}~\,$p_1>1$ ~when~ $p_1< p_2.$

Denote by $c_n$ the $n$th Gelfand number of the Sobolev embedding
(\ref{aa}). Then

$\indent c_{n}\sim n^{-\varkappa}$,\\ where\vspace{-0.2cm}
\begin{enumerate}
\item[{\rm (i)}]\ $\varkappa=\frac\mu d$\, if\, $2\le p_1\le
p_2\le \infty$\ or\,\,$1\le p_1 = p_2< 2$,\vspace{-0.2cm}

\item[{\rm (ii)}]\ $\varkappa =\frac \mu
d+\frac 1{p_1}-\frac 1{p_2}$\, if\,
$\tilde{p}<p_2<p_1\leq\infty$,\vspace{-0.2cm}

\item[{\rm (iii)}]\ $\varkappa=\frac \mu d +\frac 1{p_1}-\frac 12$\,
 if\, $1< p_1 < 2 < p_2\le\infty$\ and \,$\mu>\frac d{p_1^\prime}$,\vspace{-0.2cm}

\item[{\rm (iv)}]\ $\varkappa=\frac \mu d\cdot\frac{p_1^\prime}2$\,
 if\, $1< p_1 < 2 < p_2\le\infty$\ and \,$\mu<\frac d{p_1^\prime}$,\vspace{-0.2cm}

\item[{\rm (v)}]\ $\varkappa=\frac \mu d +\frac 1{p_1}-\frac 1{p_2}$\,
 if\, $1 < p_1 < p_2 \le 2$\ and \,$\mu>\frac d{p_1^\prime}\theta_1$,\vspace{-0.2cm}

\item[{\rm (vi)}]\ $\varkappa=\frac \mu d\cdot\frac{p_1^\prime}2$\,
 if\, $1 < p_1 < p_2 \le 2$\ and \,$\mu<\frac d{p_1^\prime}\theta_1$.
\end{enumerate}
\end{theorem}
\begin{re}\label{GBF}
In the above two theorems, the two function spaces in the embedding
(\ref{aa}) may be of different types, i.e., one is the Besov space,
and the other is the $F_{p,q}^s$-space.
\end{re}
\begin{re}
From Definition \ref{BF}, we know that an operator $f\mapsto w f$ is
an isomorphic mapping from $B_{p,q}^{s}(\mathbb{R}^d, w)$ onto
$B_{p,q}^{s}(\mathbb{R}^d)$ (similarly in the F-case). So by Remark
\ref{w12}, we get
\begin{equation*}
s_n({\rm id}, ~B_{p_1,q_1}^{s_1}(\mathbb{R}^d, v_1),~
B_{p_2,q_2}^{s_2}(\mathbb{R}^d, v_2)) ~ \sim ~ s_n({\rm id},
~B_{p_1,q_1}^{s_1}(\mathbb{R}^d, v_1/v_2),~
B_{p_2,q_2}^{s_2}(\mathbb{R}^d)),
\end{equation*}
where $s_n$ denotes either of the two quantities $c_n$ or $d_n$,
(similarly in the F-case and even in the general case emphasized in
Remark \ref{GBF}). Therefore, without loss of generality we may
assume that the target space is an unweighted space.
\end{re}
\begin{re}
For the limiting case $\delta=\alpha$, the exact order of related
n-widths may possibly depend on $q_1$ and $q_2$. In general, the
discretization method is not perfect for this case. There are
partial results on approximation numbers in \cite{MO88,Sk06} and
two-sided estimates with minor gaps in \cite{Ha97}. Some ideas from
\cite{HT05,KLSS06a} may also be helpful to further research in this
situation.
\end{re}

Now, we wish to compare the approximation, Gelfand and Kolmogorov
numbers of the Sobolev embedding (\ref{aa}). First, recall some
basic facts about approximation numbers. We define the $n$th
approximation number of $T$ by
\begin{equation}\label{an}
a_n(T)=\inf\{\|T - L\|:~ L\in \mathcal{L}(X,Y), ~\,{\rm rank} (L) <
n\},\quad n\in \mathbb{N},
\end{equation}
where ${\rm rank} (L)$ denotes the dimension of $L(X)$. We refer to
\cite{ET96,Pie78,Pin85} for detailed discussions of this concept and
further references. Let us mention that the approximation numbers
are the largest among all $s$-numbers. There exist the following
relationships:
\begin{equation}\label{acd}
 c_n(T)\le a_n(T),~~~~ d_n(T)\le a_n(T),~~~ n \in \mathbb{N}.
\end{equation}
We again assume that $1\le p_1, p_2, q_1,
q_2\le\infty,\,-\infty<s_2<s_1<\infty,\,\alpha>0,\,\delta=s_1-s_2-d(\frac
1{p_1}-\frac 1{p_2})>0,\ \delta\neq\alpha,\
\mu=\min(\alpha,\delta)$, and\ $\frac 1{\tilde{p}}=\frac\mu d+\frac
1{p_1}$. We would like to discuss, when  $a_n\sim c_n,$\ $a_n\sim
d_n,$\ or\ $c_n\sim d_n$ holds true for the Sobolev embedding
(\ref{aa}). The comparison of our results with the known results for
the approximation numbers from \cite{Sk05} shows that
\vspace{-0.2cm}
\begin{enumerate}
\item[$(i)$]\ $a_n\sim c_n$ if either

$(a)$\ $2\le p_1< p_2\le \infty$\ or,

$(b)$\ $\tilde{p}<p_2\le p_1\le \infty$\ or,

$(c)$\ $1< p_1 < p_1^\prime\le p_2\le\infty$\ and \ $\mu\neq\frac
d{p_1^\prime}$; \vspace{-0.2cm}
\item[$(ii)$]\ $a_n\sim d_n$ if either

$(a)$\ $1\le p_1< p_2\le 2$\ or,

$(b)$\ $\tilde{p}<p_2\le p_1\le \infty$\ or,

$(c)$\ $1\le p_1 < 2 < p_2\le p_1^\prime\le\infty,$\, $ p_2<
\infty$\ and\ $\mu\neq\frac d{p_2}$;\vspace{-0.2cm}
\item[$(iii)$]\ $c_n\sim d_n$ if either

$(a)$\ $\tilde{p}<p_2\le p_1\le \infty$\ or,

$(b)$\ $1< p_1 < p_1^\prime= p_2<\infty$\ and\ $\mu\neq\frac
d{p_2}$.
\end{enumerate}

\begin{re}
In the case of more general weight classes, we will show some
asymptotic estimates of the approximation, Gelfand and Kolmogorov
numbers of the corresponding embeddings in a forthcoming paper based
on ideas from T. K$\ddot{u}$hn et al. \cite{KLSS06b}.
\end{re}
\section{Sequence spaces and Kolmogorov numbers}
\subsection{Discretization of function spaces}
~~~~There are various ways to associate to a Besov space a certain
sequence space. Here we are going to use the discrete wavelet
transform, a well-developed method of discretization, see
\cite{HT05} where quasi-Banach case is included, and \cite{Tr08} for
a survey.

Let $\widetilde{\varphi}$ be an orthogonal scaling function on
$\mathbb{R}$ with compact support and of sufficiently high
regularity and let $\widetilde{\psi}$ be a corresponding wavelet.
Then the tensor product gives a scaling function $\varphi$ and
associated wavelets $\psi_1,\ldots,\psi_{2^d-1}$, all defined on
$\mathbb{R}^d$. More exactly, we suppose
$$\widetilde{\varphi}\in C^r(\mathbb{R})\quad
{\rm and \quad supp}\,\widetilde{\varphi}\subset [-N_1, N_1]$$ for
some $r\in \mathbb{N}$ and $N_1>0$. Then we have
\begin{equation}
\varphi, \psi_i \in C^r(\mathbb{R}^d) \quad{\rm and \quad
supp}\,\varphi,\, {\rm supp}\,\psi_i \subset [-N_2, N_2]^d,\quad
i=1,\ldots,2^d-1.
\end{equation}
We shall use the standard abbreviations
\begin{equation}
\varphi_{j,k}(x)=2^{jd/2}\varphi(2^jx-k) \quad{\rm and} \quad
\psi_{i,j,k}=2^{jd/2}\psi_i(2^jx-k),
\end{equation}
where $j\in\mathbb{N}_0:=\mathbb{N}\cup\{0\}$ and
$k\in\mathbb{Z}^d$.
\begin{prop}\label{Besov_des}
Let $s\in \mathbb{R}$ and $1\leq p,q\leq\infty.$ Assume $$r>\max(s,
\frac{2d}p+\frac d2-s).$$ Then a distribution $f\in \mathcal
{S}^\prime(\mathbb{R}^d)$ belongs to $B_{p,q}^{s}(w_\alpha)$ if and
only if
\begin{equation}\label{Besov_ell}
\begin{split}
\|f|B_{p,q}^{s}(w_\alpha)\|^\clubsuit&=\Big(
\sum\limits_{k\in\mathbb{Z}^d}|\langle f,\varphi_{0,k}\rangle
w_\alpha(k)|^p \Big)^{1/p}
\\
&+\sum\limits_{i=1}^{2^d-1}\Big\{ \sum\limits_{j=0}^\infty
2^{j\big(s+d(\frac 12-\frac 1p)\big)q}
\Big(\sum\limits_{k\in\mathbb{Z}^d}|\langle f,\psi_{i,j,k}\rangle
w_\alpha(2^{-j}k)|^p\Big)^{q/p}\Big\}^{1/q}<\infty.
\end{split}
\end{equation}
Moreover, $\|f|B_{p,q}^{s}(w_\alpha)\|^\clubsuit$ may be used as an
equivalent norm in $B_{p,q}^{s}(w_\alpha)$.
\end{prop}
\begin{re}
The proof of this proposition may be found in Haroske and Triebel
\cite{HT05}. One can also consult \cite{KLSS06a} for historical
remarks.
\end{re}
Let $1\leq p,q\leq\infty.$ Inspired by Proposition \ref{Besov_des} we will
work with the following weighted sequence spaces
\begin{equation}\label{ellal}
\begin{split}
\ell_q(2^{js}\ell_p(\alpha)):=\Bigg\{&
\lambda=(\lambda_{j,k})_{j,k}:~~\lambda_{j,k}\in\mathbb{C},\\
&\|\lambda|\ell_q(2^{js}\ell_p(\alpha))\|= \Big(
\sum\limits_{j=0}^\infty 2^{jsq}
\Big(\sum\limits_{k\in\mathbb{Z}^d}|\lambda_{j,k}\,w_{j,k}|^p
\Big)^{q/p}\Big)^{1/q}<\infty \Bigg\},
\end{split}
\end{equation}
(usual modification if $p=\infty$ and/or $q=\infty$), where
$w_{j,k}=w_\alpha(2^{-j}k).$ If $s=0$ we will write
$\ell_q(\ell_p(\alpha))$. In contrast to the norm defined in
(\ref{Besov_ell}), the finite summation on $i=1, 2, \ldots, 2^d-1$
is irrelevant and can be omitted.

\subsection{Kolmogorov numbers of embeddings of some sequence spaces}\label{kns}
~~~~ To begin with, we shall recall some lemmata. Lemma \ref{kn1}
follows trivially from results of Gluskin \cite{Gl83} and Edmunds
and Triebel \cite{ET96}.
\begin{lemma}\label{kn1}
Let $N\in\mathbb{N}$.
\begin{enumerate}
\item[{\rm (i)}]\ If $1\le p_1\le p_2\le 2$ and $n\le\frac{N}{4}$,\, then
$$d_{n}\left({\rm id}, \ell_{p_1}^N, \ell_{p_2}^N\right)\sim 1.$$\vspace{-0.8cm}
\item[{\rm (ii)}]\ If $1\le p_1 < 2 < p_2< \infty$ and $n\le\frac{N}{4}$,\,
then
$$d_{n}\left({\rm id}, \ell_{p_1}^N, \ell_{p_2}^N\right)\sim \min
\{1,N^{\frac 1{p_2}}n^{-\frac 1{2}}\}.$$\vspace{-0.8cm}
\item[{\rm (iii)}]\ If $2<p_1 = p_2\le \infty$ and $n\le N$,\, then
$$d_{n}\left({\rm id}, \ell_{p_1}^N, \ell_{p_2}^N\right)\sim 1.$$\vspace{-0.8cm}
\item[{\rm (iv)}]\ If $2 \le p_1 < p_2 < \infty$ and $n\le N$,\, then
$$d_{n}\left({\rm id}, \ell_{p_1}^N, \ell_{p_2}^N\right)
\sim \xi^\theta,$$ where $\xi=\min\{1,N^{\frac 1{p_2}}n^{-\frac
1{2}}\},\,\theta = \frac{1/{p_1}-1/{p_2}}{1/2-1/{p_2}}$.
\end{enumerate}
\end{lemma}

For $p_2 < p_1$ the corresponding Kolmogorov numbers can be
calculated trivially by virtue of Pietsch \cite{Pie78, Pie87}, or
Pinkus \cite[p. 203]{Pin85}.
\begin{lemma}\label{kn2} Let $1\leq p_2<p_1\leq \infty$\ and\ $n\le N$. Then
$$d_{n}\left({\rm id}, \ell_{p_1}^N, \ell_{p_2}^N\right)
=(N-n+1)^{\frac 1{p_2}-\frac 1{p_1}}.$$
\end{lemma}

The following assertion is a simple corollary of Lemma \ref{kn1}.
And the proof is the same as that of Lemma 10 in Skrzypczak
\cite{Sk05}.
\begin{lemma}\label{kn3}  Suppose $1\le p_1 < 2 < p_2< \infty$ and $N\in\mathbb{N}$.
 Then
there is a positive constant $C$ independent of $N$ and $n$ such
that
\begin{equation}\label{klp2p}
d_{n}\left({\rm id},\ell_{p_1}^N,\ell_{p_2}^N\right)\leq C\,
\begin{cases}
1 & {\rm if}\ n\leq N^{\frac 2{p_2}},\\
N^{\frac 1{p_2}}n^{-\frac 1{2}}\quad & {\rm if}\ N^{\frac
2{p_2}}<n\leq
N,\\
0 & {\rm if}\ n>N.
\end{cases}
\end{equation}
\end{lemma}

\begin{prop}\label{kn5}  Suppose $1\le p_1 < 2 < p_2< \infty$
 and $\delta\neq\alpha$. We set
\begin{equation}\label{kp2pp}
\varkappa\, =\, \begin{cases}
\frac \mu d +\frac 12-\frac 1{p_2} &{\rm if}\,\, \mu>\frac d{p_2},\\
\frac \mu d\cdot\frac{p_2}2 &{\rm if}\,\, \mu<\frac d{p_2},
\end{cases}
\end{equation}where $\mu=\min(\alpha,\delta).$ Then
\begin{equation}d_{n}\left({\rm id}, \ell_{q_1}(2^{j\delta}\ell_{p_1}(\alpha)),
\ell_{q_2}(\ell_{p_2})\right) \sim n^{-\varkappa}.\end{equation}
\end{prop}
\noindent{\bf Proof.} By Lemma \ref{kn1} and Lemma \ref{kn3}, the
proof of the proposition can be finished in the same manner as in
the proof of Prop.\,11 in \cite{Sk05} with the important complement
given in \cite{SV09}. The only change is that
$t=\min(p_1^\prime,p_2)$ is replaced by $p_2$ in our proof. \qed
\begin{prop}\label{kn6}  Suppose $2 \le p_1 < p_2 < \infty$
 and $\delta\neq\alpha$. We set
\begin{equation}\label{kp2pp}
\varkappa\, =\, \left\{
\begin{array}{ll}
\frac \mu d +\frac 1{p_1}-\frac 1{p_2} &{\rm if}\,\, \mu>\frac d{p_2}\theta,\\
\frac \mu d\cdot\frac{p_2}2 &{\rm if}\,\, \mu<\frac d{p_2}\theta,
\end{array}\right.
\end{equation}where $\mu=\min(\alpha,\delta),\,\theta = \frac{1/{p_1}-1/{p_2}}{1/2-1/{p_2}}.$ Then
\begin{equation}d_{n}\left({\rm id}, \ell_{q_1}(2^{j\delta}\ell_{p_1}(\alpha)),
\ell_{q_2}(\ell_{p_2})\right) \sim n^{-\varkappa}.\end{equation}
\end{prop}
\noindent{\bf Proof.} {\tt Step 1}.\quad Preparations. We denote
$$\Lambda:=\{\lambda=(\lambda_{j,k}):
\,\quad\lambda_{j,k}\in\mathbb{C},\,\quad j\in
\mathbb{N}_0,\,k\in\mathbb{Z}^d\},$$ and set
$$ B_1=\ell_{q_1}(2^{j\delta}\ell_{p_1}(\alpha))\,\quad {\rm
and}\,\quad  B_2= \ell_{q_2}(\ell_{p_2}).$$
Let\,$I_{j,i}\subset\mathbb{N}_0\times\mathbb{Z}^d$ be such that
\begin{equation}\label{Ij0}
I_{j,0}:=\{(j,k):\, |k|\leq 2^j\},\quad j\in\mathbb{N}_0,
\end{equation}
\begin{equation}\label{Iji}
I_{j,i}:=\{(j,k):\, 2^{j+i-1}<|k|\leq 2^{j+i}\}, \quad i\in
\mathbb{N},\quad j\in\mathbb{N}_0.
\end{equation} Besides, let
$P_{j,i}:\Lambda\mapsto\Lambda$ be the canonical projection onto the
coordinates in $I_{j,i}$; i.e., for $\lambda\in\Lambda$, we set
\begin{equation*}
(P_{j,i}\lambda)_{u,v}:= \left\{
\begin{array}{ll}
\lambda_{u,v}\quad & (u,v)\in I_{j,i},\\
0\quad & {\rm otherwise},
\end{array}\right. \quad u\in
\mathbb{N}_0,\quad v\in\mathbb{Z}^d,\quad i\ge 0.
\end{equation*}
Then
\begin{equation}\label{rankpiece}M_{j,i}:=|I_{j,i}|\sim2^{d(j+i)},\end{equation}
\begin{equation}{\rm id}_\Lambda=
\sum\limits_{j=0}^\infty\sum\limits_{i=0}^\infty
P_{j,i},\end{equation}
\begin{equation}w_{j,k}=w_\alpha(2^{-j}k)\sim2^{\alpha i}
\quad {\rm if}\quad (j,k)\in I_{j,i},\ \ i\ge 0.\end{equation} Due
to simple monotonicity arguments and explicit properties of the
Kolmogorov numbers we have
\begin{equation}
\begin{array}{ll}\label{dndiscr}
d_n(P_{j,i},B_1,B_2)&\leq  \frac 1{\inf_{k\in
I_{j,i}}w_{j,k}}2^{-j\delta}d_n({\rm id}, \ell_{p_1}^{M_{j,i}}, \ell_{p_2}^{M_{j,i}})\\
&\leq c2^{-j\delta-i\alpha}d_n({\rm id}, \ell_{p_1}^{M_{j,i}},
\ell_{p_2}^{M_{j,i}}).
\end{array}
\end{equation}
{\tt Step 2}.\quad The operator ideal comes into play. Under the
assumption $2 \le p_1 < p_2 < \infty$, it is easy to prove that
$0<\theta\le 1.$ To shorten notations we shall put $\tau=\frac
{p_2}\theta,\,h=\frac 2\theta,\,$ and\,$\frac 1s=\frac 1\gamma+\frac
1h$\, for any $s>0.$\, In terms of (\ref{idealddef}) and
(\ref{dndiscr}), we have
\begin{equation}\label{idealdiscr}
L_{s,\infty}^{(d)}(P_{j,i})\leq
c2^{-j\delta-i\alpha}L_{s,\infty}^{(d)}({\rm id},
\ell_{p_1}^{M_{j,i}}, \ell_{p_2}^{M_{j,i}}).
\end{equation}
The known asymptotic behavior of the Kolmogorov numbers $d_n({\rm
id},\ell_{p_1}^{N},\ell_{p_2}^{N})$, cf. Lemma \ref{kn1} (iv), and
(\ref{rankpiece}) yield that
\begin{equation}\label{idealseh}
L_{h,\infty}^{(d)}({\rm id}, \ell_{p_1}^{M_{j,i}},
\ell_{p_2}^{M_{j,i}})\leq
C2^{d(j+i)/\tau},\indent\indent\indent\indent\indent\quad
\end{equation}
\begin{equation}\label{idealsh}
L_{s,\infty}^{(d)}({\rm id}, \ell_{p_1}^{M_{j,i}},
\ell_{p_2}^{M_{j,i}})\leq C2^{d(j+i)(\frac 1 \tau+\frac
1\gamma)}\quad\quad\ \,\, {\rm if} \quad \frac 1s>\frac 1h.
\end{equation}
{\tt Step 3}.\quad The estimate of $d_n({\rm id}, B_1, B_2)$ from
above in the first case $\mu>\frac d\tau.$\, For any given
$M\in\mathbb{N}_0$,\, we put
\begin{equation}\label{PQ}
P:=\sum\limits_{m=0}^M\sum\limits_{j+i=m}P_{j,i}\quad\quad{\rm
and}\quad\quad
Q:=\sum\limits_{m=M+1}^\infty\sum\limits_{j+i=m}P_{j,i}.
\end{equation}
{\tt Substep 3.1.}\quad Estimate of $d_n(P, B_1, B_2)$. Let $\frac
1s>\frac 1h$. Then in view of (\ref{idealdinq}),\,(\ref{idealdiscr})
and (\ref{idealsh}), we have
\begin{equation}
\begin{aligned}\label{idealPrho}
L_{s,\infty}^{(d)}(P)^\rho&\leq
\sum\limits_{m=0}^M\sum\limits_{j+i=m}L_{s,\infty}^{(d)}(P_{j,i})^\rho\\
&\leq c_1
\sum\limits_{m=0}^M\sum\limits_{j+i=m}2^{-\rho(j\delta+i\alpha)}2^{\rho
md(\frac 1\tau+\frac 1\gamma)}\\
&\leq c_2 \sum\limits_{m=0}^M 2^{\rho md(\frac 1\tau+\frac
1\gamma-\frac \mu d)}.
\end{aligned}
\end{equation}
In the last inequality we used our assumption $\delta\neq \alpha$.
We choose $\gamma$ such that $d(\frac 1\tau+\frac 1\gamma)-\mu>0.$
Then (\ref{idealPrho}) yields
\begin{equation}
\label{idealP} L_{s,\infty}^{(d)}(P)\leq c~ 2^{dM(\frac 1\tau+\frac
1\gamma-\frac \mu d)}.
\end{equation}
Using (\ref{idealddef}) and (\ref{idealP}), we get
\begin{equation}
\label{d_2dMP} d_{2^{dM}}(P, B_1, B_2)\leq c_32^{dM(\frac
1\tau-\frac 1h-\frac \mu d)}.
\end{equation}
We put $n=2^{Md}.$ Then
\begin{equation}
\label{d_nP} d_n(P, B_1, B_2)\leq c_3n^{\frac 1\tau-\frac 1h-\frac
\mu d}=c_3n^{-(\frac \mu d+\frac 1{p_1}-\frac 1{p_2})}.
\end{equation}
{\tt Substep 3.2.}\quad Estimate of $d_n(Q, B_1, B_2)$. In a similar
way to (\ref{idealPrho}), we obtain by $\delta\neq\alpha$ and
(\ref{idealseh}) that
\begin{equation}
\label{idealQrho} L_{h,\infty}^{(d)}(Q)^\rho\leq c_1
\sum\limits_{m=M+1}^\infty2^{\rho md(\frac 1\tau-\frac \mu d)}.
\end{equation}
Since $\mu>\frac d\tau,$\, we have
\begin{equation}
\label{idealQ} L_{h,\infty}^{(d)}(Q)\leq c 2^{dM(\frac 1\tau-\frac
\mu d)}.
\end{equation}
$\indent$ By virtue of (\ref{idealddef}), we have
\begin{equation}
\label{d_2dMQ} d_{2^{dM}}(Q, B_1, B_2)\leq c_22^{dM(\frac
1\tau-\frac 1h-\frac \mu d)}.
\end{equation}
Take $n=2^{Md}.$ Then
\begin{equation}
\label{d_nQ} d_n(Q, B_1, B_2)\leq c_2n^{\frac 1\tau-\frac 1h-\frac
\mu d}=c_2n^{-(\frac \mu d+\frac 1{p_1}-\frac 1{p_2})}.
\end{equation}
{\tt Substep 3.3.}\quad Under the assumption $\mu>\frac d\tau,$\,
the estimate from above follows from (\ref{d_nP}) and (\ref{d_nQ}),
by means of the inequality below,
\begin{equation}
\label{d_nid} d_{2n}({\rm id}, B_1, B_2)\leq d_n(P, B_1, B_2)+d_n(Q,
B_1, B_2).
\end{equation}
{\tt Step 4}.\quad The estimate of $d_n({\rm id}, B_1, B_2)$ from
above in the second case $\mu<\frac d\tau.$ Inspired by \cite{SV09},
we use the following division
\begin{equation}
\label{divi3} {\rm
id}=\sum\limits_{m=0}^{M_1}\sum\limits_{j+i=m}P_{j,i}+
\sum\limits_{m=M_1+1}^{M_2}\sum\limits_{j+i=m}P_{j,i}+
\sum\limits_{m=M_2+1}^\infty\sum\limits_{j+i=m}P_{j,i},
\end{equation}
where $M_1, M_2\in\mathbb{N}$\,\,and \,$M_1<M_2$, which will be
determined later on. Using the subadditivity of $s$-numbers, we have
\begin{equation}
\label{dn_divi3} d_{n^\prime}({\rm id}, B_1, B_2)\leq
\triangle_1+\triangle_2+\triangle_3,
\end{equation}
where
\begin{equation*}
\begin{array}{ll}
\triangle_1=\sum\limits_{m=0}^{M_1}
\sum\limits_{j+i=m}d_{n_{j,i}}(P_{j,i}),\quad\quad &
\triangle_2=\sum\limits_{m=M_1+1}^{M_2}
\sum\limits_{j+i=m}d_{n_{j,i}}(P_{j,i}), \\
\triangle_3=\sum\limits_{m=M_2+1}^\infty\sum\limits_{j+i=m}\|P_{j,i}\|,\quad\quad
&n^\prime-1=\sum\limits_{m=0}^{M_2}\sum\limits_{j+i=m}(n_{j,i}-1).
\end{array}
\end{equation*}
Note that for $\triangle_3$, we have $j+i > M_2$ in the sum, and we
take $n_{j,i}=1.$ Now let $n\in\mathbb{N}$ be given. We take
\begin{equation*}
M_1=\left[\frac {\log_2n}d-\frac{\log_2\log_2n}d\right]\quad\quad
{\rm and}\quad \quad M_2=\left[\frac \tau
h\cdot\frac{\log_2n}d\right],
\end{equation*}
where $[a]$ denotes the largest integer smaller than
$a\in\mathbb{R}$ and $\log_2n$ is a dyadic logarithm of $n$. Then
\begin{equation*}
\begin{array}{ll}
\triangle_3&=\sum\limits_{m=M_2+1}^\infty\sum\limits_{j+i=m}\|P_{j,i}\|
\leq
c_1\sum\limits_{m=M_2+1}^\infty\sum\limits_{j+i=m}2^{-j\delta}2^{-i\alpha}
\\ &\leq c_2\sum\limits_{m=M_2+1}^\infty2^{-m\mu}\leq
c_32^{-M_2\mu}\leq c_3n^{-\varkappa}.
\end{array}
\end{equation*}
Next, we choose proper $n_{j,i}$ for estimating $\triangle_1$ and
$\triangle_2$. If $i+j\leq M_1,$ we take $n_{j,i}=M_{j,i}+1$ such
that $d_{n_{j,i}}(P_{j,i})=0$ and $\triangle_1=0$. And we obtain
$$\sum\limits_{m=0}^{M_1}\sum\limits_{j+i=m}n_{j,i}\leq
 c_1\sum\limits_{m=0}^{M_1}(m+1)2^{md}\leq c_2M_1 \cdot 2^{M_1d}\leq c_3 n.$$
Now we give the crucial choice of $n_{j,i}$ for the second sum
$\triangle_2.$ We take
$$n_{j,i}=[n^{1-\varepsilon}\cdot2^{iz_1}\cdot2^{jz_2}],$$
where $\varepsilon, z_1, z_2$ are positive real numbers such that
$$\alpha+\frac{z_1}h<\frac d\tau,\quad
0<\frac{z_1-z_2}h<\delta-\alpha\quad {\rm and}\quad
\frac{z_1\tau}{hd}=\varepsilon\quad {\rm if}\, \delta>\alpha,$$ or
$$\delta+\frac{z_2}h<\frac d\tau,\quad
0<\frac{z_2-z_1}h<\alpha-\delta\quad {\rm and}\quad
\frac{z_2\tau}{hd}=\varepsilon\quad {\rm if}\, \delta<\alpha.$$ Note
that the relation, $0<\varepsilon<1,$ holds obviously. Then
$$\sum\limits_{m=M_1+1}^{M_2}
\sum\limits_{j+i=m}n_{j,i}\leq
c_1n^{1-\varepsilon}\sum\limits_{m=M_1+1}^{M_2}2^{m\cdot\max(z_1,z_2)}
\leq c_2n^{1-\varepsilon}\cdot
n^{\frac\tau{hd}\max(z_1,z_2)}=c_2n,$$ and, in terms of
(\ref{idealseh}),
\begin{equation*}
\begin{aligned}
\sum\limits_{m=M_1+1}^{M_2}\sum\limits_{j+i=m}d_{n_{j,i}}(P_{j,i})
&\leq
c_1\sum\limits_{m=M_1+1}^{M_2}\sum\limits_{j+i=m}2^{-j\delta-i\alpha}
2^{(i+j)d/\tau}[n^{1-\varepsilon}\cdot2^{iz_1}\cdot2^{jz_2}]^{-\frac 1h}\\
&\leq c_2n^{-\frac 1h(1-\varepsilon)}\sum\limits_{m=M_1+1}^{M_2}
2^{md/\tau}2^{-m\cdot\min(\alpha+\frac{z_1}h,\delta+\frac{z_2}h)}\\
&\leq c_3n^{-\frac 1h(1-\varepsilon)}n^{\frac
1h}n^{-\frac\tau{hd}\min(\alpha+\frac{z_1}h,\delta+\frac{z_2}h)}\\
&= c_3n^{\frac\varepsilon
h-\frac\tau{hd}\min(\alpha+\frac{z_1}h,\delta+\frac{z_2}h)}\\
&= c_3n^{-\frac{\tau\mu}{hd}}= c_3n^{-\varkappa}.
\end{aligned}
\end{equation*}
Hence the estimate from above in the second case is finished.
\\{\tt Step 5}.\quad The estimate of $d_n({\rm id}, B_1, B_2)$ from
below. Consider the following diagram
\begin{equation}
\begin{CD}
\ell_{p_1}^{M_{j,i}} @>S_{j,i}>>
\ell_{q_1}(2^{j\delta}\ell_{p_1}(\alpha))
\\
@VV{\rm id_1}V @VV{\rm id}V\\
\ell_{p_2}^{M_{j,i}} @<T_{j,i}<< \ell_{q_2}(\ell_{p_2})
\end{CD}
\end{equation}
Here,
\begin{equation*}
\begin{array}{rl}
(S_{j,i}\eta)_{u,v}:=&\left\{
\begin{array}{ll}\eta_{\varphi(u,v)}\,
&{\rm if}\quad (u,v)\in I_{j,i}, \\ 0 \, &{\rm otherwise,}
\end{array}\right.\vspace{0.1cm}\\
(T_{j,i}\lambda)_{\varphi(u,v)}:=&\lambda_{u,v},\quad\quad\quad
(u,v)\in I_{j,i},\end{array}
\end{equation*}
and $\varphi$ denotes a bijection of $I_{j,i}$ onto $\{1, \ldots,
M_{j,i}\},\, j\in\mathbb{N}_0,\, i\in\mathbb{N}_0$; cf. (\ref{Ij0})
and (\ref{Iji}). Observe that
\begin{equation*}
\begin{array}{ll}
S_{j,i}\in \mathcal {L}\left(\ell_{p_1}^{M_{j,i}},\,
\ell_{q_1}(2^{j\delta}\ell_{p_1}(\alpha))\right) \quad  &{\rm
and}\quad \|S_{j,i}\|=2^{j\delta+i\alpha},\vspace{0.1cm}\\
T_{j,i}\in  \mathcal {L}\left(\ell_{q_2}(\ell_{p_2}),\,
\ell_{p_2}^{M_{j,i}}\right) \quad  &{\rm and}\quad \|T_{j,i}\|=1.
\end{array}
\end{equation*}
Hence we obtain
\begin{equation}\label{diagST}
d_n(\,{\rm id}_1)\leq\|S_{j,i}\|\,\|T_{j,i}\|\,d_n(\,{\rm id}).
\end{equation}

(i)\, Let $\frac d\tau<\delta\leq\alpha.$  We consider
$N:=M_{j,0}=|I_{j,0}|\sim2^{d j},\,j\ge\frac 2d.$ Then
$$\|S_{j,0}\|\leq C2^{j\delta}\quad{\rm and}\quad\|T_{j,0}\|=1.$$
Put $m=\frac N4\sim2^{jd-2}.$ And for sufficiently large $N$ we have
$m\ge N^{\frac 2{p_2}}$ since $p_2>2.$ Consequently,
\begin{equation*}
d_m({\rm id}_1,\,\ell_{p_1}^N,\,\ell_{p_2}^N)\sim\left(N^{\frac
1{p_2}}m^{-\frac 12}\right)^\theta\sim2^{\theta(jd-2)(\frac
1{p_2}-\frac 12)}\sim2^{(jd-2)(\frac 1{p_2}-\frac 1{p_1})}.
\end{equation*}
Using (\ref{diagST}), we obtain
\begin{equation*}
d_{2^{jd-2}}({\rm id})\ge C_12^{-j\delta}2^{(jd-2)(\frac
1{p_2}-\frac 1{p_1})}\ge C_22^{(jd-2)(\frac 1{p_2}-\frac
1{p_1}-\frac\delta d)}.
\end{equation*}
Then the monotonicity of the Kolmogorov numbers implies that for any
$ n\in \mathbb{N}$
\begin{equation}
d_n({\rm id})\ge C_3n^{-(\frac\delta d+\frac 1{p_1}-\frac 1{p_2})}.
\end{equation}

(ii)\, Let $\frac d\tau<\alpha<\delta.$ We consider
$N:=M_{0,i}=|I_{0,i}|\sim2^{d i},\,i\ge\frac 2d.$ Then
$$\|S_{0,i}\|\leq C2^{i\alpha}\quad{\rm and}\quad\|T_{0,i}\|=1.$$
Also put $m=\frac N4\sim2^{di-2}.$ Hence we have similarly for any $
n\in \mathbb{N}$
\begin{equation}
d_n({\rm id})\ge Cn^{-(\frac\alpha d+\frac 1{p_1}-\frac 1{p_2})}.
\end{equation}

(iii)\, Let $\delta\leq\frac d\tau\ {\rm and}\ \delta<\alpha$. We
select the same $N,\,S$, and\ $T$ as in point (i) and take
$m=\left[N^{\frac 2{p_2}}\right]\leq\frac N4$ for sufficiently large
$N.$ Then $N^{\frac 1{p_2}}m^{-\frac 12}\sim 1.$ Hence by Lemma
\ref{kn1} and (\ref{diagST}) we obtain
\begin{equation*}
d_m({\rm id})\ge C2^{-j\delta}=C2^{-jd\frac
2{p_2}\frac{p_2\delta}{2d}},
\end{equation*}
and then for any $ n\in \mathbb{N}$
\begin{equation}
d_n({\rm id})\ge Cn^{-\frac{p_2\delta}{2d}}.
\end{equation}

(iv)\, Let $\alpha\leq\frac d\tau\ {\rm and}\ \alpha\leq\delta$. We
select the same $N,\,S$, and\ $T$ as in point (ii) and take
$m=\left[N^{\frac 2{p_2}}\right]$ in the same way as in point (iii).
Then analogously
\begin{equation*}
d_m({\rm id})\ge C2^{-i\alpha}=C2^{-di\frac
2{p_2}\frac{p_2\alpha}{2d}},
\end{equation*}
and in consequence, for any $ n\in \mathbb{N}$
\begin{equation}
d_n({\rm id})\ge Cn^{-\frac{p_2\alpha}{2d}}.
\end{equation}
The proof of the proposition is now complete. \qed
\begin{prop}\label{kn7}  Suppose $1\le p_1\le p_2\le
2$\,\,or\,\,$2<p_1 = p_2\le \infty$\,and $\delta\neq\alpha.$ We set
\begin{equation}\varkappa
=\frac \mu d\,,\quad\quad
where\quad\mu=\min(\alpha,\delta).\end{equation} Then
\begin{equation}d_{n}\left({\rm id}, \ell_{q_1}(2^{j\delta}\ell_{p_1}(\alpha)),
\ell_{q_2}(\ell_{p_2})\right) \sim n^{-\varkappa}.\end{equation}
\end{prop}

In view of Lemma \ref{kn1}, the proof of this proposition follows
exactly as in the proof of Prop.\,13 in \cite{Sk05}.
\begin{prop}\label{kn8}  Suppose $1\le p_1, p_2\le\infty,\,\,
\frac 1{\tilde{p}}=\frac{\min(\alpha,\delta)}d+\frac 1{p_1},$\,and
$\delta\neq\alpha.$ Assume $\tilde{p}<p_2<p_1\leq\infty$, and set
\begin{equation}\varkappa =\frac \mu d+\frac 1{p_1}-\frac
1{p_2}\,,\quad\quad where\quad\mu=\min(\alpha,\delta).\end{equation}
Then
\begin{equation}d_{n}\left({\rm id}, \ell_{q_1}(2^{j\delta}\ell_{p_1}(\alpha)),
\ell_{q_2}(\ell_{p_2})\right) \sim n^{-\varkappa}.\end{equation}
\end{prop}

By Lemma \ref{kn2}, the proof of this proposition can be finished in
the same manner as in the proof of Prop.\,15 in \cite{Sk05}.
\section{Proofs of the main results}
\subsection{Proof of Theorem \ref{kn}}\label{proofkn}
~~~\, Based on Prop.~\ref{Besov_des}, we now transfer the results of
Subsection \ref{kns} for weighted sequence spaces to weighted
function spaces.

First, for the embeddings given by (\ref{emB_wei}), i.e., the Besov
case, the assertions follow from Proposition~\ref{Besov_des} and
\ref{kn5}-\ref{kn7}, or
Proposition~\ref{kn8}, respectively.

For the general case, we estimate from above, by virtue of the
multiplicativity property of Kolmogorov numbers, and the elementary
embeddings below
\begin{equation*}
A_{p_1,q_1}^{s_1}(\mathbb{R}^d,w_\alpha) \hookrightarrow
B_{p_1,\infty}^{s_1}(\mathbb{R}^d,w_\alpha) \hookrightarrow
B_{p_2,1}^{s_2}(\mathbb{R}^d) \hookrightarrow
A_{p_2,q_2}^{s_2}(\mathbb{R}^d).
\end{equation*}
For the estimate from below we can consider the following embeddings
\begin{equation*}
B_{p_1,1}^{s_1}(\mathbb{R}^d,w_\alpha) \hookrightarrow
A_{p_1,q_1}^{s_1}(\mathbb{R}^d,w_\alpha) \hookrightarrow
A_{p_2,q_2}^{s_2}(\mathbb{R}^d) \hookrightarrow
B_{p_2,\infty}^{s_2}(\mathbb{R}^d).
\end{equation*}\qed
\subsection{Proof of Theorem \ref{gn}}
~~~~ We turn our attention to Gelfand numbers. First, we collect
some necessary information on $c_n({\rm id}, \ell_{p_1}^N,
\ell_{p_2}^N)$, cf. \cite{Gl83, Pie78, Vy08}, (\ref{dualc*d}) and
(\ref{duald*c}).
\begin{lemma}\label{gn1} Let $N\in\mathbb{N}$.
\begin{enumerate}
\item[{\rm (i)}]\ If $2\le p_1\le
p_2\le\infty$ and $n\le\frac{N}{4}$,\, then
$$c_{n}\left({\rm id}, \ell_{p_1}^N, \ell_{p_2}^N\right)\sim 1.$$
\vspace{-0.8cm}
\item[{\rm (ii)}]\ If $1< p_1 < 2 < p_2\le \infty$ and $n\le\frac{N}{4}$,\,
then
$$c_{n}\left({\rm id}, \ell_{p_1}^N, \ell_{p_2}^N\right)\sim \min\{1,N^{1-\frac 1{p_1}}n^{-\frac 1{2}}\}.$$
\vspace{-0.8cm}
\item[{\rm (iii)}]\ If $1\le p_1 = p_2< 2$\ and $n\le N$,\, then
$$c_{n}\left({\rm id}, \ell_{p_1}^N, \ell_{p_2}^N\right)\sim 1.$$
\vspace{-0.8cm}
\item[{\rm (iv)}]\ If $1 < p_1 < p_2 \le 2$\ and $n\le N$,\, then
$$c_{n}\left({\rm id}, \ell_{p_1}^N, \ell_{p_2}^N\right)
\sim \xi^{\theta_1},$$ where $\xi=\min\{1,N^{1-\frac
1{p_1}}n^{-\frac 1{2}}\},\,\theta_1 =
\frac{1/{p_1}-1/{p_2}}{1/{p_1}-1/2}$.
\end{enumerate}
\end{lemma}
The proof of this lemma follows by (\ref{dualc*d}), (\ref{duald*c})
and Lemma \ref{kn1}.
\begin{lemma}\label{gn2}
Let $1\leq p_2<p_1\leq \infty$\ and\ $n\le N$. Then
$$c_{n}\left({\rm id}, \ell_{p_1}^N, \ell_{p_2}^N\right)
=(N-n+1)^{\frac 1{p_2}-\frac 1{p_1}}.$$
\end{lemma}
The proof of this lemma follows literally Pietsch
\cite{Pie78,Pie87}, and also Pinkus \cite{Pin85}.

Now we are ready to prove Theorem \ref{gn}.  We can deal with the
proof in a similar way to the one for Theorem \ref{kn}, so we give
the sketch here. For example, in the case (v), which is
corresponding to the case (v) in Theorem \ref{kn}, the proof can be
based on Prop.~\ref{kn6}. The changes begin with (\ref{dndiscr}),
where $d_n$ is substituted by $c_n$. And the others go on trivially.
\begin{re}
For the quasi-Banach case, $0<p<1 ~{\rm or}~0<q<1$, the problem of
these two quantities of the embeddings given by (\ref{aa}) becomes
more complicated. Indeed, Lemma \ref{kn1} and Lemma \ref{kn2} can
not be completely generalized to the quasi-Banach setting $0<p_1,
p_2\le \infty$. And the duality between these two quantities is not
valid. Fortunately, some recent results provided by Foucart et al.
\cite{FPRU} and Vyb\'iral \cite{Vy08} are effective to a certain
extent. This situation will be discussed in a forthcoming paper.
\end{re}

\section*{Acknowledgments}

~~~~ The authors wish to thank the anonymous referees and Professor
K. Ritter for their excellent comments, remarks and suggestions
which greatly helped us to improve this paper. The authors are also
extremely grateful to Ant\'onio M. Caetano, Fanglun Huang, Thomas
K$\ddot{\rm u}$hn, Erich Novak and Leszek Skrzypczak for their
direction and help on this work.

This research was partially supported by the Natural Science
Foundation of China (Grant No. 10671019, No. 11171137) and Anhui Provincial
Natural Science Foundation (No. 090416230).

\end{document}